\documentclass[11pt]{amsart}
\usepackage{amsfonts,amssymb,amsthm,amsmath,amsxtra,amscd,verbatim,eucal}
\usepackage[all]{xy}
\usepackage[dvips]{graphics}

\theoremstyle{plain}
\newtheorem{thm}{Theorem}[section]
\newtheorem{lem}[thm]{Lemma}
\newtheorem{prop}[thm]{Proposition}

\newtheorem{conj}[thm]{Conjecture}

\theoremstyle{definition}
\newtheorem{defi}[thm]{Definition}

\theoremstyle{remark}
\newtheorem{rmk}[thm]{Remark}

\newtheorem*{ack}{Acknowledgements}

\numberwithin{equation}{section}

\def\C{{\mathbb C}}

\def\R{{\mathbb R}}

\def\P{{\mathbb P}}

\def\B{\mathcal{B}}

\def\O{\mathcal{O}}

\def\mm{\mathfrak{m}}

\def\a{\alpha}

\def\d{\delta}

\def\ff{\psi}

\def\ep{\epsilon}

\def\l{\lambda}

\def\p{\pi}
\def\r{\rho}
\def\s{\sigma}
\def\t{\tau}

\def\G{\Gamma}

\def\S{\Sigma}
\def\Om{\Omega}

\def\.{\cdot}
\def\^{\widehat}
\def\~{\widetilde}
\def\o{\circ}
\def\ov{\overline}

\def\de{\partial}

\def\({\left(}
\def\){\right)}

\DeclareMathOperator{\rk} {rk}

\DeclareMathOperator{\NE} {NE}
\DeclareMathOperator{\CNE} {\ov{\NE}}

\DeclareMathOperator{\Nef} {Nef}

\DeclareMathOperator{\Pic} {Pic}

\DeclareMathOperator{\Nul} {Nul}

\DeclareMathOperator{\mult} {mult}

\DeclareMathOperator{\Div} {Div}

\DeclareMathOperator{\diag} {diag}

\DeclareMathOperator{\Pos} {Pos}

\begin{document}

\title{Negative curves on very general blow-ups of $\P^2$}

\thanks{Research partially supported by
the University of Michigan Rackham Research Grant and the MIUR of
the Italian Government, National Research
Project ``Geometry on Algebraic Varieties" (Cofin 2002).}

\author{Tommaso de Fernex}
\address{Department of Mathematics, University of Michigan,
East Hall, 525 East University Avenue, Ann Arbor, MI 48109-1109, USA}
\email{{\tt defernex@umich.edu}}

\begin{abstract}
This note contains new evidence to a conjecture, related to the Nagata
conjecture and the Segre-Harbourne-Gimigliano-Hirschowitz conjecture,
on the cone of effective curves of blow-ups of $\P^2$ at
very general points.
\end{abstract}


\maketitle

\section*{Introduction}

The solution to the problem of determining the dimension of every linear
system of curves in $\P^2$ with assigned multiplicities
at general points is predicted by the equivalent conjectures
of Segre, Harbourne, Gimigliano and Hirschowitz
\cite{Seg61,Har86,Gim87,Hir89}, hereafter ``SHGH conjecture".
In this paper we are interested in the following weaker form of the conjecture:

\begin{conj}\label{conj1}
Suppose that $C$ is an integral curve with negative
self-intersection on the blow-up $Y$ of $\P^2$ at a set
of points in very general position. Then $C$ is a $(-1)$-curve of $Y$
(that is, a smooth rational curve with self-intersection $-1$).
\end{conj}

If one also includes the conjecture that $H^1(Y,\O_Y(D))=0$
for every effective nef divisor $D$ on any such $Y$, then this
becomes one of the formulations of the SHGH conjecture \cite{Har2}.
There is an interesting relationship between
Conjecture~\ref{conj1}, which of course contains Nagata's conjecture \cite{Nag59},
and the symplectic packing problem in dimension four \cite{MP,Bir}.
In the first section of this paper, we will review how the above
conjecture characterizes the geometry of the cone of effective curves of $Y$.

It is not difficult to see that the conjecture is satisfied
by all rational curves (see Proposition~\ref{rational} below).
There are several results nowadays, giving evidence to the SHGH conjecture,
that are valid under suitable assumptions on the multiplicities
assigned to the centers of the blow-up. These include
\cite{AH,Har3,CM98,CM00,Mig,Yan}.
In the same spirit, we prove the following result.

\begin{thm}\label{main-intro}
Conjecture~\ref{conj1} is satisfied by
all curves whose image on $\P^2$ is a curve
with a singularity of multiplicity 2
at one of the centers of the blow-up.
\end{thm}

The proof of this theorem is based on a local study
of dynamic self-intersection, very much
inspired to the methods in \cite{EL} and \cite{Xu}.
The novelty here with respect to \cite{EL,Xu} is to consider deformations
with two-parameters families: it is indeed by computing the Kodaira-Spencer
map from different directions of deformation that we produce a
linear system on $C$  giving an isomorphism to $\P^1$ and
showing that $C^2 = -1$. A similar result is
proven to hold for an arbitrary smooth projective surface in place of $\P^2$
(see Theorem~\ref{main} below).

An analogous result on
linear systems on smooth projective surfaces was given
in \cite[Theorem~4.1]{dVL} under the assumption
that the ``specialty'' of the system
(namely, the gap between its dimension and its {\it expected dimension})
increases by imposing a double point. We remark that
this hypothesis does not translate well to the context of
determining negative curves on the blow-up.
To the best of our knowledge, there are no other results in which
the assumptions only involve one of the multiplicities.

The precise notion of ``very general position"
adopted in this paper is given in Definition~\ref{defi} below.
Throughout this paper we work over the complex numbers.

\begin{ack}
The author would like to thank L.~Ein, B.~Harbourne,
A.~Laface, R.~Lazarsfeld, R.~Miranda and S.~Yang for useful discussions.
\end{ack}

\section{Geometry of the cone of curves}\label{geometry}

In this section we discuss the implication that
Conjecture~\ref{conj1} has on the
geometry of the cone of effective curves of the blow-up of $\P^2$.
This section is mostly of an expository nature, as the material
here contained is probably well known to the specialist.

We start by fixing some notation.
Let $X$ be a smooth projective surface.
Let
$$
N(X) := (\Pic(X)/\equiv)\otimes \R
$$
and $\r := \rk \Pic(X)$.
By the Hodge index theorem, we can identify $N(X)$ with
$\R^\r$ in such a way that the intersection product is given by the
matrix $\diag(1,-1,\dots,-1)$. For a divisor $D$ on $X$
we denote by $D^{< 0}$ (resp. $D^{\le 0}$, $D^{\ge 0}$, $D^{\perp}$)
the subset of $N(X)$ defined by $[D]\.x < 0$
(resp. $[D]\.x \le 0$, $[D]\.x \ge 0$, $[D]\.x = 0$). Let
$$
\CNE(X) \subset N(X)
$$
be the closure of the cone spanned by the classes of
effective curves on $X$, and
$$
\Nef(X) \subset N(X)
$$
be the closure of the
cone spanned by the classes of ample divisors on $X$. By Kleiman's
criterion, these two cones are put in duality by the intersection
product in $N(X)$. Fix an ample class $H$ on $X$, let
$$
\Pos(X) := \{ \a \in N(X) \mid \a^2 \ge 0,\; \a\.H \ge 0 \},
$$
and let $\Nul(X)$ denote the boundary of $\Pos(X)$.
Then $\Nul(X)$ is supported by the quadratic
equation $x_1^2 = x_2^2 + \dots + x_{\r}^2$. Note that
$$
\Nef(X) \subseteq \Pos(X) \subseteq \CNE(X).
$$
It follows by a result of Campana and Peternell~\cite{CP} that
$\de \CNE(X)$ is supported by $\Nul(X)$ and countably many hyperplanes,
so we can write
$$
\de \CNE(X) = \B_1 \sqcup \B_2,
$$
where $\B_1$ is the closure (in $\de \CNE(X)$) of
the union of the facets of $\de \CNE(X)$, and
$\B_2$ is supported by $\Nul(X)$.

The following is a more precise formulation of Conjecture~\ref{conj1}.

\begin{conj}\label{NE(Y)}
Let $Y$ be the blow-up of $\P^2$ at a set of $r$
points in very general position. Then, writing $\de \CNE(Y) = \B_1 \sqcup \B_2$
as above (setting $X = Y$), the extremal rays of $\B_1$ that do not lie on $\Nul(Y)$
are spanned by classes of $(-1)$-curves, and we have
$$
\B_1 = \de \CNE(Y) \cap K_Y^{\le 0}
\quad\text{and}\quad \B_2 = \de \CNE(Y) \cap K_Y^{> 0}.
$$
In particular,
$$
\CNE(Y) \cap K_Y^{\ge 0} = \Pos(Y) \cap K_Y^{\ge 0}.
$$
\end{conj}

In the following, we review the equivalence of the two conjectures.
We begin with a general consideration on the clustering
of extremal rays. As at the beginning of the section,
let $X$ be a smooth projective surface.
We consider the metric on the set of rays in $N(X)$ given by the
angular distance: for any two rays $R_1$ and $R_2$ in $N(X)$, we
set $d(R_1,R_2)$ to be the angle between them.
For an extremal ray $R$ of $\CNE(X)$, we set
$$
d(R) = \inf \{d(R,R') \mid \text{$R'$ is an extremal ray of $\CNE(X)$
different from $R$} \}.
$$

\begin{lem}\label{cluster}
If $R$ is an extremal ray of $\CNE(X)$, then $d(R,\Nul(X)) \le d(R)$;
in particular,
every Cauchy sequence of extremal rays of $\CNE(X)$ converges to a
ray on $\Nul(X)$. Moreover, if $\r \ge 4$ and $F$ is a
facet of $\CNE(X)$ such that $\de F \cap \Nul(X) \ne \emptyset$,
then $F$ has infinitely many extremal rays forming a Cauchy sequence
converging to $\de F \cap \Nul(X)$.
\end{lem}

\begin{proof}
We can assume that $R \not \subset \Nul(X)$.
By \cite[Lemma~II.4.2]{Kol}, $R$ is spanned by the class of a curve $C$
with $C^2 < 0$, and every other extremal ray of $\CNE(X)$
is contained in the half space $C^{\ge 0}$ \cite[Exercise~1.4.33(ii)]{Laz}.
Thus the inequality $d(R,\Nul(X)) \le d(R)$ will follow once
we show that the orthogonal projection (in
the euclidian metric) of $R$ onto $C^{\perp}$ is contained in $\Pos(X)$.
This can be easily checked as follows.
Let $\langle \;,\; \rangle$ denote the standard
inner product in $\R^{\r}$, and let $(\;\.\;)$ denote the intersection
product defined by the diagonal matrix $\diag(1,-1,\dots,-1)$.
Let $h = (1,0,\dots,0) \in \R^{\r}$, fix a vector
$c \in \R^{\r}$ such that $(c\.h) > 0$ and $(c\.c) < 0$, and consider
the orthogonal projection
$$
\p : \R^{\r} \to \{x \in \R^{\r} \mid (c\.x) = 0 \}.
$$
It is an exercise to check that
$$
\frac{(\p(c)\.\p(c))}{(c\.c)} = \frac{(c\.c)^2}{\langle c,c \rangle^2} - 1
\quad\text{and}\quad
(\p(c)\.h) = (c\.h) - \frac{(c\.c)}{\langle c,c \rangle}\.\langle c,h \rangle.
$$
Since $(c\.c)^2 \le \langle c,c \rangle^2$, the first equation gives
$(\p(c)\.\p(c)) \ge 0$. The second equation gives $(\p(c)\.h) > 0$.
Therefore the projection of $R$ onto $C^\perp$ in contained in $\Pos(X)$.

To prove the second part
of the lemma, we note that $F$ is supported by the hyperplane $H$ that is
tangent to $\Nul(X)$ along the ray $\de F \cap \Nul(X)$; this
follows by the convexity of $\CNE(X)$
and the inclusion $\Pos(X) \subseteq \CNE(X)$. If $F$ has
only finitely many rays, then we find another facet $F'$
containing $\de F \cap \Nul(X)$ that, for the same reason, is also supported by $H$.
Since this is impossible, $F$ must have infinitely many extremal rays,
and these cluster to $\de F \cap \Nul(X)$
by what proved in the first part of the lemma.
\end{proof}

Now we can go back to the question on the equivalence between the
two conjectures. Recall that we are denoting by
$Y$ the blow-up of $\P^2$ at a set of $r$
points in very general position.
Clearly both conjectures are certainly true for $r \le 2$,
so we can assume that $r \ge 3$, hence $\rk \Pic(Y) \ge 4$.
Since one direction is obvious, let us assume that the only integral curves
$C \subset Y$ with $C^2 < 0$ are the $(-1)$-curves.
Write $\de \CNE(Y) = \B_1 \sqcup \B_2$ as described above.
The extremal rays of $\B_1$ not lying on $\Nul(Y)$,
which are also extremal rays of $\CNE(Y)$,
are spanned by classes of $(-1)$-curves.
In particular, by adjunction formula and Lemma~\ref{cluster}, we have
$\B_1 \subset K_Y^{\le 0}$. On the other hand, the Cone Theorem
implies that $\B_2 \subset K_Y^{\ge 0}$.
In fact, observing that $\de\CNE(Y) \cap K_Y^{< 0} \ne \emptyset$
and that $\B_2$ is open in $\de\CNE(Y)$, we
actually have $\B_2 \subset K_Y^{> 0}$.
Therefore we conclude that
$$
\B_1 = \de \CNE(Y) \cap K_Y^{\le 0}\quad\text{and}\quad
\B_2 = \de \CNE(Y) \cap K_Y^{> 0}.
$$
The last assertion follows by continuity.

\begin{rmk}
As it is well known, an interest in Conjecture~\ref{NE(Y)} comes from the
search for examples of irrational Seshadri constants.
Indeed, if true, the conjecture would allow us to construct
such examples. This is easy to see: let $Y$ be the blow-up of $\P^2$ at $r \ge 9$
points in very general position, and denote by $H$ the pull-back
on $Y$ of the hyperplane class of $\P^2$ and by $E$ the exceptional
divisor of $Y \to \P^2$. Then, for a very general point $q \in Y$
and for every rational $0 < a < 1/\sqrt{r}$, we would have
$\ep_q(H - aE) = \sqrt{1 - ra^2}$.
Note that, for most (rational) values of $a$, this
is an irrational number.
\end{rmk}

\section{Negative curves on $\P^2$ blown up at
very general points}\label{negative}

The necessity to assume very general position in the conjecture is clear
in the case of more than 9 points. Indeed the blow-up of $\P^2$
at 9 points will generally carry infinitely $(-1)$-curves \cite[Theorem~4a]{Nag60},
and certainly one cannot allow to blow up points lying on these curves.
In this paper we adopt the following notion of very general position.

\begin{defi}\label{defi}
We say that a set $\{x_1,\dots,x_r\}$
of $r \ge 1$ distinct points on a smooth projective surface $X$
is {\it in very general position}
if the following condition is fulfilled.
For every integral curve $C \subset X$, the pair
$(C, (x_1,\dots,x_r))$ belongs to an irreducible algebraic family
$$
\{ (C_t,(x_{1,t},\dots,x_{r,t})) \mid t\in T \},
$$
where $C_t \subset X$ is an integral curve and
$(x_{1,t},\dots,x_{r,t}) \in X^r$ for every $t \in T$,
satisfying the following properties:
\begin{enumerate}
\item[(i)]
$p_g(C_t) = p_g(C)$ and
$\mult_{x_{i,t}}C_t = \mult_{x_i}C$ for every $t \in T$ and $i=1,\dots,r$;
\item[(ii)]
the morphism $\ff : T \to X^r$ given by $\ff(t) = (x_{1,t},\dots,x_{r,t})$
is an isomorphism to an open subset of $X^r$.
\end{enumerate}
We say that a point $x \in X$ is a {\it very general point}
if the set $\{x\}$ is in very general position.
\end{defi}

\begin{rmk}\label{rmk-very}
For every family of integral curves on $X$, the set of $r$-ples
$(x_1,\dots,x_r) \in X^r$
for which the condition of the definition is not satisfied
by some curve of the family
is contained in a proper closed subvariety of $X^r$.
We conclude that the complement of the locus in $X^r$
of points corresponding to subsets of $X$ in
very general position is a countable union of proper
closed subvarieties.
Any subset of a set of points in very general position is
a set of points in very general position.
Given a set of points in very general position,
the image of any of the points of the set on the blowing-up
of the surface at the residual points of the set is a very general point.
\end{rmk}

The blow-up of $\P^2$ at at most 8 points in general position
is a Del Pezzo surface. This surface contains
finitely many $(-1)$-curves, and these are the only integral curves
with negative self-intersection.
The picture in the case of $9$ points is well known too
\cite[Proposition~12]{Nag60};
a proof of the following statement will also be given below.

\begin{prop}\label{r=9}
Every integral curve with negative
self-intersection on the blow-up of $\P^2$
at a set of 9 points in very general position is a $(-1)$-curve.
\end{prop}

Very little is known on the cone of effective curves
of the blow-up of $\P^2$ at more than 9 points.
If we restrict our attention to rational curves, then
it is easy to get the desired statement.
The following proposition, which is well known to the
experts, can be viewed as the first step
with respect to another formulation of the SHGH conjecture,
due to Harbourne \cite{Har04}, which says that
if $C$ is any irreducible and reduced curve on
the blow-up of $\P^2$ at very general points, then
$C^2 \ge p_a(C) - 1$, where $p_a(C)$ is the arithmetic genus of $C$.

\begin{prop}\label{rational}
Let $Y$ be the blow-up of $\P^2$ at a set of points in very general position,
and assume that $C$ is a integral rational curve on $Y$
with $C^2 < 0$. Then $C$ is a $(-1)$-curve.
\end{prop}

\begin{proof}
To prove that $C$ is a $(-1)$-curve, it suffices
to show that $K_Y\.C < 0$, as the conclusion then follows by adjunction.
Let $B \subset \P^2$ be the image of $C$. We can assume that $B$ is
a curve, otherwise $C$ would be exceptional, hence a $(-1)$-curve.
Let $p_1,\dots,p_r \in \P^2$ be the centers of the blow-up, and
let $m_i = \mult_{p_i}B$. We can assume that $m_1 > 0$.

For short, let $p = (p_1,\dots,p_r)$,
By the definition of very general position, the pair
$(B, p)$ belongs to an irreducible algebraic family
$$
\{ (B_t, p_t) \mid t\in T \},
$$
where $p_t = (p_{1,t},\dots,p_{r,t}) \in (\P^2)^r$,
$B_t \subset X$ is an integral rational curve with
$\mult_{p_{i,t}}B_t = m_i$ for every $t \in T$, and
the morphism $\ff : T \to (\P^2)^r$ given by $\ff(t) = p_t$
is an isomorphism to an open subset of $(\P^2)^r$.

Fix $r$ points $q_1,\dots,q_r$ in general position
on a smooth cubic $\G \subset \P^2$.
Let $Z \subset (\P^2)^r$ be a smooth, irreducible curve passing
through $p$ and $q=(q_1,\dots,q_r)$, and consider the open set
$U = Z \cap\ff(T)$ of $Z$. Note that $U$ is not empty, and that
$q$ is in the closure of $U$.
The family $\{(B_t,p_t)\mid t \in \ff^{-1}(U)\}$ determines
an effective divisor
$$
\B \subset \P^2 \times U,
$$
whose restriction to $\P^2 \times \{p_t\}$ (for every $t \in \ff^{-1}(U)$)
is the divisor $B_t$.
Viewing $\B$ as a scheme,
we take its flat closure $\ov\B$ inside $\P^2 \times Z$,
and let $B_0$ be the restriction of $\B$ to $\P^2 \times \{q\}$.
Then $B_0$ is an effective divisor on $\P^2$
(see for instance \cite[Example~III.9.8.5]{Ha}).
Since $\mult_{p_{i,t}}B_t = m_i$ for every $p_t \in U$,
we have $\mult_{q_i}B_0 \ge m_i$ by the semi-continuity of the multiplicity.
Moreover, since $B_t$ is a rational curve for every $t \in \ff^{-1}(U)$,
so is every irreducible component
of $B_0$, hence $\G$ is not contained in the support of $B_0$. Thus
\begin{equation}\label{B}
-K_{\P^2}\.B = \G\.B_0 \ge \sum m_i.
\end{equation}

If this inequality is strict, then we have
$$
K_Y\.C = K_{\P^2}\.B + \sum m_i < 0.
$$
Therefore, to conclude the proof, it is enough to
show that the inequality in~(\ref{B}) is strict.
Suppose this is not the case. Then $\mult_{q_i}B_0 = m_i$ and
$\O_{\P^2}(B_0)|_{\G} = \O_{\G}(\sum m_iq_i)$.
Moreover, by the way we chose the $q_i$,
we can fix a different deformation in which $q_1$
is replaced by another general point $q_1'$ of $\G$ while the other
$q_i$ are kept the same, obtaining in this way another curve $B_0'$.
Since the previous arguments also apply to $B_0'$, we have
$\mult_{q_1'}(B_0') = m_1$ and $\mult_{q_i}(B_0') = m_i$ for $i \ge 2$, and
moreover
$$
\O_{\G}\(m_1q_1 + \sum_{i=2}^rm_iq_i\) = \O_{\P^2}(B_0)|_{\G} =
\O_{\P^2}(B_0')|_{\G} = \O_{\G}\(m_1q_1' + \sum_{i=2}^rm_iq_i\).
$$
This implies that $\O_{\G}\(m_1(q_1 - q_1')\) = \O_{\G}$.
But $\G$ is an elliptic curve, and after fixing $q_1$ as zero,
we know that there are finitely many $m_1$-torsion points on $\G$.
A contradiction.
\end{proof}

\begin{proof}[Proof of Proposition~\ref{r=9}.]
Let $Y$ be the blow-up of $\P^2$ at a set of 9 points in
very general position, and suppose that $C \subset Y$
is an integral curve with $C^2 < 0$.
Since $-K_Y$ is nef, we have $K_Y\.C \le 0$.
By adjunction, we conclude that $C$ is rational, hence it is
a $(-1)$-curve by Proposition~\ref{rational}.
\end{proof}

The following is the main result of this paper.
Choosing $S = \P^2$ in the statement gives the
theorem stated in the introduction.

\begin{thm}\label{main}
Let $S$ be a smooth projective surface, and
let $f : Y \to S$ be the blowing up of $S$ at a set
$\S$ of points in very general position.
Let $C \subset Y$ be an integral curve with negative
self-intersection, and assume that $f(C)$
is a curve with a singularity of multiplicity 2
at one of the points of $\S$. Then $C$ is a $(-1)$-curve of $Y$.
\end{thm}

\begin{proof}
Let $p \in \S$ be the double point of $f(C)$ whose existence we are assuming,
and let $X$ be the blow-up of $S$ at $\S \setminus \{p\}$.
Then $f$ factors through
the blow-up $g : Y \to X$ of the image $x \in X$ of $p$.
Let $D = g(C)$. Note that $\mult_x D = 2$.
Since $x$ is a very general point of $X$,
there exists an irreducible algebraic family of integral curves
with marked points
$$
\{ (D_t,x_t) \mid t\in T \} \subset X \times X
$$
with $D = D_{t^*}$ and $x = x_{t^*}$ for some $t^* \in T$,
such that the morphism $\ff : T \to X$ given by $\ff(t) = x_t$ is an isomorphism
to an open subset of $X$ and, moreover, $\mult_{x_t}D_t = 2$ for all $t \in T$.
The Kodaira-Spencer map induced by any 1-dimensional degeneration $t \to t^*$
inside $T$ defines a section of the normal
bundle $N := N_{D/X}$ of $D$ in $X$.
Bearing in mind that we are dealing with a family of irreducible
and reduced curves with marked singularities, these sections are non-zero
whenever the degeneration $t\to t^*$ is performed along
a curve of $T$ that is smooth at $t^*$.

As in the proof of \cite[Lemma~1.1]{EL}, we reduce to a local computation in
some open set $\Om$ in $\C^2$. We
fix local coordinates $u = (u_1,u_2)$ in $\Om$.
Let $f = f(u)$ be the holomorphic function locally defining $D$.
Let us start considering the case in which $p$ is an ordinary node.
Writing $f$ as a power series centered at $(0,0)$,
we can assume that the coordinates are chosen so that
\begin{equation}\label{f}
f(u) = u_1u_2 + \text{(higher degree terms)}.
\end{equation}
We can reduce to the case when $T$ is a small disk in $\C^2$,
with $t^* = (0,0)$, and fix coordinates $t = (t_1,t_2)$ in $T$
such that $t_i = u_i\ff$. The total space of the deformation
is defined in $\Om \times T$ by a power series $F = F(u,t)$.
The deformation determines a Kodaira-Spencer map
$$
\r : T_{t^*}(T) \to H^0(D,N),
$$
which is non-trivial by our previous assumptions.
This map is locally defined by
\begin{equation}\label{tau}
\r\left(\l_1\frac{\de\;}{\de t_1} + \l_2\frac{\de\;}{\de t_2}\right) =
\left(\l_1\frac{\de F}{\de t_1}(u,0) +
\l_2\frac{\de F}{\de t_2}(u,0)\right)\Big{|}_C =:\t_{\l}.
\end{equation}
In view of the linearity of this map, the sections
$\t_{\l}$ fill up a non-trivial linear subspace in $H^0(D,N)$
as $\l$ varies in $\C^2$.
Mimicking~\cite{EL}, we consider the function
$$
\Phi(u,t) := F(u + x(t), t),
$$
where $x(t)=(x_1(t),x_2(t))$ are the coordinate
of the marked point $x_t$ of $C_t$.
Note that $\Phi \in (u_1,u_2)^2$ for all $t$.
We expand $\Phi$ as a power series in $(t_1,t_2)$.
The coefficients of the two terms of degree 1 are equal to
\begin{equation}\label{Phi}
\frac{\de \Phi}{\de t_i} =
\frac{\de f}{\de u_1}(u) \.\frac{\de x_1}{\de t_i}(0) +
\frac{\de f}{\de u_2}(u) \.\frac{\de x_2}{\de t_i}(0) +
\frac{\de F}{\de t_i}(u,0), \quad i=1,2.
\end{equation}
These are functions of $u$, and both are contained in
$(u_1,u_2)^2$. Note that $\de x_j/\de t_i(0) = \d_{ij}$.
Then, by combining~(\ref{f}),~(\ref{tau}) and~(\ref{Phi}), we see that
\begin{equation}\label{m}
(\l_1u_2 + \l_2u_1)|_C + \t_{\l} \in \mm_x^2,
\end{equation}
where $\mm_x$ is the maximal ideal of $x$ in $D$.
Note that $(\l_1u_2 + \l_2u_1)|_C \in \mm_x$. We conclude that
$\t_{\l} \in H^0(D,N\otimes \mm_x)$.
Then $\t_{\l}$ gets pulled back by $g|_C : C \to D$
to a section $\s_{\l}$ of $(g|_C)^*N$ that
vanishes at the two pre-images of $x$.
After suitably denoting these two pre-images by $y_1$ and $y_2$,
we actually get sections
$$
\s_{(1,0)}\in H^0(C,(g|_C)^*N\otimes \mm_{y_1}^2 \otimes \mm_{y_2})
\quad\text{and}\quad
\s_{(0,1)}\in H^0(C,(g|_C)^*N\otimes \mm_{y_1} \otimes \mm_{y_2}^2)
$$
when $\l \in \{(1,0),(0,1)\}$.
This implies that $\deg(\Div(\s_{\l})) \ge 3$ for every $\l \ne 0$.
Since $D^2 < 4$, this yields
$\deg (\Div(\s_{\l})) = 3$, that is $D^2 = 3$, and thus $C^2 = -1$.
In fact, we observe that the linear system $|\Div(\s_{\l})|$ contains
a pencil parameterized by $\l$. This pencil
has base points at $y_1$ and $y_2$ and movable part of degree 1,
which defines an isomorphism to $\P^1$.

It remains to discuss the case when $x$ is not an ordinary node of $D$.
We now explain why this case cannot occur.
Using analogous notation as in the previous discussion,
we get the following local equation of $D$:
$$
f = u_1^2 + \text{(higher degree terms)}.
$$
Arguing as before, we see this time that
$\t_{(0,1)} \in H^0(C,N\otimes \mm_x^2)$.
This implies that $\deg(\Div(\s_{(0,1)})) \ge 4$, which is
impossible.
\end{proof}

\providecommand{\bysame}{\leavevmode \hbox \o3em
{\hrulefill}\thinspace}


\begin{thebibliography}{AAA}

\bibitem[AH]{AH}
J. Alexander and A. Hirschowitz, Polynomial interpolation in several variables,
\emph{J. Algebraic Geom.} \textbf{4} (1995), 201--222.

\bibitem[Bi]{Bir}
P. Biran, Symplecting packing in dimension 4,
{\it Geom. Funct. Anal.} {\bf 7} (1997), 420--437.

\bibitem[CP]{CP}
F. Campana and T. Peternell, Algebraicity of the ample
cone of projective varieties, \emph{J. Reine Angew. Math.}
\textbf{407} (1990), 160--166.

\bibitem[CM1]{CM98}
C. Ciliberto and R. Miranda, Degenerations of planar linear systems,
\emph{J. Reine Angew. Math.} \textbf{501} (1998), 191--220.

\bibitem[CM2]{CM00}
C. Ciliberto and R. Miranda, Linear systems of plane curves
with base points of equal multiplicity,
\emph{Trans. Amer. Math. Soc.} \textbf{352} (2000), 4037--4050.

\bibitem[dVL]{dVL}
C. de Volder and A. Laface, Linear systems on generic K3 surfaces,
to appear in {\it Bull. Belg. Math. Soc.}, {\tt math.AG/0309073}.

\bibitem[EL]{EL}
L. Ein and R. Lazarsfeld, Seshadri constants on smooth surfaces,
\emph{Ast\'{e}risque} (1993), no. 218, 177--186, Journ\'{e}es de
G\'{e}om\'{e}trie Alg\'{e}brique d'Orsay (Orsay, 1992).

\bibitem[Gi]{Gim87}
A. Gimigliano, On linear systems of plane curves,
Ph.D. thesis, Queen's University (1987).

\bibitem[Ha1]{Har86}
B. Harbourne, The geometry of rational surfaces and Hilbert
functions of points in the plane,
\emph{Can. Math. Soc. Conf. Proc.} \textbf{6} (1986), 95--111.

\bibitem[Ha2]{Har2}
B. Harbourne, Points in good position in $\P^2$,
in {\it Zero-Dimensional Schemes},
Proceedings of the International Conference,
Ravello, 1992 (F. Orecchia and L. Chiantini Eds.),
de Guyter, pp. 213--229, 1994.

\bibitem[Ha3]{Har3}
B. Harbourne, Rational surfaces with $K^2 > 0$,
{\it Proc. Amer. Math. Soc.} {\bf 124} (1996), 727--733.

\bibitem[Ha4]{Har04}
B. Harbourne, The (unexpected) importance of knowing $\a$,
preprint.

\bibitem[Har]{Ha}
R. Hartshorne, \emph{Algebraic Geometry},
Graduate Texts in Mathematics, No. 52, Springer-Verlag, New York, 1977.


\bibitem[Hi]{Hir89}
A. Hirschowitz, Une conjecture pour la cohomologie
des diviseurs sur les surfaces rationelles g\'en\'eriques,
\emph{J. Reine Angew. Math.} \textbf{397} (1989), 208--213.

\bibitem[Ko]{Kol}
J. Koll\'ar, {\em Rational Curves on Algebraic Varieties},
Ergeb. Math. Grenzgeb. (3) 32, Springer-Verlag, Berlin, 1996.

\bibitem[Laz]{Laz}
R. Lazarsfeld, {\it Positivity in Algebraic Geometry, I},
Ergeb. Math. Grenzgeb. (3) 48, Springer-Verlag, Berlin, 2004.

\bibitem[MP]{MP}
D. McDuff and L. Polterovich, Symplecting packings and
algebraic geometry, With an appendix by Y. Karshon,
{\it Invent. Math.} {\bf 115} (1994), 405--434.

\bibitem[Mi]{Mig}
R. Migon, Syst\`eme de courbes planes \`a singularit\'es impos\'ees:
le cas des multiplicit\'es inf\'erieures ou \'egales \`a quatre,
\emph{J. Pure Appl. Algebra} \textbf{151} (2000), 173--195.

\bibitem[Na1]{Nag59}
M. Nagata, On the 14-th problem of Hilbert,
\emph{Amer. J. Math.} \textbf{81} (1959), 766--772.

\bibitem[Na2]{Nag60}
M. Nagata, On rational surfaces, II,
\emph{Mem. Coll. Sci. Univ. Kyoto, Ser.~A},
Vol. XXXIII, Math. No. 2 (1960), 271--293.

\bibitem[Se]{Seg61}
B. Segre, Alcune questioni su
insiemi finiti di punti in Geometria Algebrica,
Atti del Convegno Internazionale di Geometria Algebrica, Torino (1961).

\bibitem[Ya]{Yan}
S. Yang, Linear series in $\P^2$ with base points of bounded multiplicity,
preprint, {\tt math.AG/0406591}.

\bibitem[Xu]{Xu}
G. Xu, Curves in $\P^2$ and symplectic packings,
\emph{Math. Ann.} (1994), 609--613.


\end{thebibliography}
\end{document}